\documentclass[twoside,12pt,reqno]{amsart}

\usepackage{upref,amsxtra,amssymb,amscd}
\usepackage{varioref}
\usepackage{verbatim}
\usepackage{epsfig}

\usepackage{epic,eepic,eucal}
\usepackage{enumerate}

\def\B{          \mathcal B}

\def\A{          \mathcal A}

\def\a{         \alpha}

\newcommand{\NN}{{\mathbb N}}
\newcommand{\RR}{{\mathbb R}}
\newcommand{\TT}{{\mathbb T}}

\newcommand{\ZZ}{{\mathbb Z}}

\newcommand{\QQ}{{\mathbb Q}}

\def\carre{ \hfill $\Box$    }

\newtheorem{theo}{\sc Theorem}[section]
\newtheorem{prop}[theo]{\sc Proposition}
\newtheorem{lemm}[theo]{\sc Lemma}
\newtheorem{coro}[theo]{\sc Corollary}

\theoremstyle{definition}

\newtheorem{defi}[theo]{\sc Definition}

\theoremstyle{remark}

\newtheorem{rema}[theo]{\sc Remark}

\numberwithin{equation}{section}

\begin{document}
\title{Two remarks on the Shrinking
  Target Property}
\author[Bassam Fayad]{Bassam Fayad}
\address{Bassam  Fayad, LAGA, Universit\'e Paris 13, Villetaneuse, }
\email{fayadb@math.univ-paris13.fr}

\begin{abstract}

  We show that a translation $T_\a$ of the torus $\TT^d$
has  the {\it monotone
    shrinking target property} if and only if  the vector 
    $\a$ is of constant type. 

Then, using reparametrizations of linear flows, we show that there
exist area preserving real analytic maps of the three torus that are mixing of all orders and do not enjoy the 
  {\it monotone
    shrinking target property}.
\end{abstract}

\maketitle

\section{Introduction}

Let $(T,M,\mu)$ be an ergodic dynamical system. For any measurable set $A \subset M$, such that  $\mu(A) > 0$, there exist for almost every $x \in M$ infinitely
many integers $n$ such that $T^n(x) \in A$. If $A$ is  called the {\it target}, we can say that the trajectory
of almost every point $x$ hits  the {target} infinitely often. In an attempt to refine this recurrence property due to ergodicity, we can consider instead of a fixed set $A$ a sequence of sets $A_n$ such that  $\mu(A_n) \rightarrow 0$ and investigate whether 
the events $T^n(x) \in A_n$ occur infinitely many times for a.e. $x \in M$. Of course, $\mu(A_n)$ cannot converge to zero too fast, since the trivial implication of the Borel-Cantelli Lemma asserts that if $\sum \mu(A_n) < \infty$, then for almost every $x \in A$, $T^n(x) \in A_n$ occurs only finitely many times.

A general abstract definition exists:     
\begin{defi} \label{bc} A sequence of measurable sets ${\A}= {\lbrace  A_n \rbrace}_{n \in \NN }$  such that 
\begin{eqnarray} \label{div} \sum_{n=0}^{\infty} \mu (A_n) =  \infty, \end{eqnarray}
is called a {\it Borel-Cantelli (BC) sequence} for $T$ if the sets 
$$ H(x, {\A})  =  \lbrace n \in \NN \ / \ T^n(x) \in A_n \rbrace $$ 
are infinite for almost every point $x$ in $M$. In other words, $\A$ is BC if and only if 
$$\mu \left(\lim \sup T^{-n}A_n\right)=1,$$
where $\lim \sup T^{-n}A_n  = \bigcap_{n_0 \in \NN} \bigcup_{n \geq n_0} T^{-n}A_n.$
\end{defi} 

The points $x \in \lim \sup T^{-n}A_n $ are said to hit the shrinking target infinitely often (we refer to \cite{kleinbock-chernov} for more definitions and a nice historical account on the use of the above formalism in dynamics). 

This definition was motivated  by the classical theory of Diophantine approximations. From the point of view of dynamics, we see by different versions of the Borel-Cantelli Lemma that if the sets $T^{-n}A_n$ satisfy some independence conditions  then  $\A$ is BC under the assumption (\ref{div}). Checking the above Borel-Cantelli property is actually a measurement of the stochasticity of the system $(T,M,\mu)$ (see \S  \ref{blabla} below).

Given a dynamical system $(T,M,\mu)$, obtaining {\it Dynamical Borel-Cantelli lemmas} consist in proving that for a sequence of targets in some particular familly of sets, condition (\ref{div}) implies that the sequence is BC. Indeed, condition (\ref{div}) by itself is not enough to guarantee that $ H(x, {\A})$ is infinite for a.e. $x \in M$ since 
 for any measure preserving
system $(T,M,\mu)$, there exists a sequence $\A$ satisfying (\ref{div}) that is not  a BC sequence (the sequence can even be taken nested, that is $A_{n+1} \subset A_n$, for all $n\geq 0$). So in order to characterize dynamical systems using BC sequences we need to be more specific and put some restrictions on the nature  of the subsets that form the target. In the context of ergodic continuous maps  $T$ on a
metric space $M$, it is natural to consider the following definition: 

\begin{defi}[STP] \label{tttSTP} We say that  $(T,M,\mu)$ has the
{\it Shrinking Target Property} (STP) if for any $x_0 \in M$,  any sequence of
balls centered at $x_0$ that satisfies (\ref{div}) is BC for
$T$. \end{defi}


We can restrict even more the choice of the target to sequences of
decreasing centered balls, $\B =  {\lbrace  B (x_0,r_n) \rbrace}_{n \in \NN }$ with   
$0 \leq r_{n+1} \leq r_n$ for $n \in \NN$, and this would
 yield a second definition of the STP, called the monotone shrinking target
property or MSTP:

\begin{defi}[MSTP] \label{tttSTP2} We say that T has the
{\it Monotone Shrinking Target Property} if for any $x_0 \in M$, any sequence of decreasing centered balls around $x_0$ that satisfies (\ref{div}) is BC for $T$. \end{defi} 

\begin{rema} Notice that if the radii of the balls in $\B$ are decreasing, then 
the set of points $x \in M$ such that $H(x,\B)$ is infinite is invariant by $T$ (since $T^nx \in B(x_0,r_n)$ implies  $T^{n-1}(Tx) \in B(x_0,r_{n-1})$). Hence, for an ergodic $T$, the set of points that hit infinitely often a monotone shrinking target is either of measure 0 or is of measure 1. Hence, to check the MSTP for  an ergodic  system $(T,M,\mu)$, it is enough  to prove that for any $x_0$ there exists $\eta>0$ such that for any decreasing sequence $r_n$ such that $\B ={\lbrace B(x_0,r_n) \rbrace}_{n \in \NN}$ satisfies (\ref{div}) we have that 
$$\mu \left( \bigcup_{n \geq 0} T^{-n} B(x_0,r_n) \right) \geq \eta.$$
\label{eta}
\end{rema}

\section{Toral translations}


\subsection{Definitions and Notations} $ \ $

\noindent $\bullet \ $ For $d \geq 1$, we consider the $d$-dimensional torus  $\TT^d = \RR^d / \ZZ^d$, and for $\a \in \RR^d$ we denote by $T_\a$ the translation of vector $\a$ on $\TT^d$. Any translation preserves the Haar measure on $\TT^d$ that we denote by $\mu^{(d)}$, or simply $\mu$ when there is no possible ambiguity.

\noindent $\bullet \ $ We denote the sup norm on $\RR^d$ by
$$|\a| = \sup (|\a_1|,\cdots,|\a_d|).$$
$\bullet \ $  For $\a \in \RR$, we denote its closest distance to the integers by
$$\|\a\| = \inf_{p \in \ZZ} |\a-p|.$$
$\bullet \ $ For $\a \in \RR^d$, define 
$${\|\a\|}_{\ZZ} = \sup_{i} \|\a_i\|.$$
$\bullet \ $ For $\a \in \RR^d$ and $k \in \ZZ^d$, denote by $(k,\a)$ the scalar $\sum_{i=1}^{d} k_i\a_i$. For $Q \in \ZZ$, let $Q\a = (Q\a_1,\cdots,Q\a_d)$.

\vspace{0.1cm} 

We recall the definition of Diophantine vectors in $\RR^d$
\begin{defi}  For $d \geq 1$ and $\sigma >0$, we denote by $\Omega^d(\sigma)$  the set of  Diophantine vectors in $\RR^d$ of exponent $\sigma d$, that is
$$\Omega^d(\sigma) = \lbrace \a \in \RR^d \ / \ \exists C>0 \ {\rm such \ that \ } \forall k \in \ZZ^d \setminus \lbrace 0 \rbrace \quad
  \|(k,\a)\| \geq C {{|k|}^{-(1+\sigma)d}}   \rbrace.$$ 
The set $\Omega_d(0)$ is said to be the set of {\it constant type vectors} in $\RR^d$.
\end{defi}

There is a similar definition based on simultaneous approximations
\begin{defi} For $d \geq 1$ and $\sigma>0$ define
$$\Omega^d(\sigma) = \lbrace \a \in \RR^d \ / \ \exists C>0 \ {\rm such \ that \ } \forall Q \in \ZZ \setminus \lbrace 0 \rbrace \quad  {\|Q\a\|}_{\ZZ}  \geq C {|Q|}^{- {(1+\sigma) \over d} }  \rbrace.$$ 
\end{defi}

The Khintchine transference  principle, Cf. for example \cite[Chapter V]{cassels}, asserts that 
$$\Omega_d \left({\sigma \over (d-1)\sigma +d} \right) \subset \Omega^d(\sigma) \subset \Omega_d(d\sigma).$$
In the particular case of constant type vectors we get the identity
$$\Omega_d(0)=\Omega^d(0).$$
In this note, we will only use the latter definition of  constant type vectors, because simultaneous approximations are naturally involved in the study of the Shrinking target properties.

\subsection{STP vs MSTP} 
 
Since an additional restriction on the target is imposed in the definition of the MSTP we have that if a system $(T,M,\mu)$ has the STP, then it necessarily has the MSTP but the converse is false as illustrated by the following: 

\begin{theo} \label{theorem1} Let $\a \in \RR^d$. Consider the translation $T_\a$ of vector $\a$  on $\TT^d$, and let $\mu$ be the Haar measure on $\TT^d$. Then we have


\begin{itemize}
\item $(T_{\a},\TT^d,\mu)$ has the MSTP if and only if the vector $\a$ is of constant type.
\end{itemize}
\begin{itemize} 
\item $\forall \a \in \RR^d$, $(T_\a,\TT^d,\mu)$ does not have the STP.
\end{itemize}
Moreover, if $\a \in \RR^d$ is not a rational vector, there exists a sequence $r_n \geq 0$ satisfying $\sum  r_n^d = \infty$ such that, for any $x_0 \in \TT^d$,  the sequence $\B = {\lbrace B(x_0,r_n)\rbrace}_{n \in \NN}$ satisfies $\lim \sup T_\a^{-n} B(x_0,r_n) = \emptyset$, that is, the set of points 
that hit the target infinitely often is empty. 
\end{theo}

\begin{rema} In $\RR^d$ (or $\TT^d$), a sequence of balls $B(x_n,r_n)$ satisfies (\ref{div}) for the Lebesgue (or Haar) measure if and only if $$\sum_{n=0}^{\infty} r_n^d = \infty.$$
\end{rema}

\begin{rema} If $\a \in \RR^d$ is rational, the set of points that do not hit the target infinitely often cannot be empty. It is easy to see however that for any $x_0 \in \TT^d$ and for any nonempty subset $\chi$ of the finite orbit of $x_0$ by $T_\a$, we can find a sequence of centered balls around $x_0$ satisfying (\ref{div}) such that the set of points 
that hit the target infinitely often is reduced to $\chi$. \end{rema}

\subsection{} We will start with the proof of the second assertion which is more elementary. Let  $\a \in \RR^d$ and assume that one of the coordinates of $\a$, say $\a_1$ is not rational. 

Let $k_p \in \NN$, $p \geq 1$, be an increasing sequence of integers such that the intervals  $[-k_p \a_1 - 4^{-p}, -k_p \a_1 + 4^{-p}]$ are disjoint.  

Let $q_n$, $n \geq 1$, be a sequence of integers such that 
$$\|q_n\a_1\| \leq e^{-n}$$
and define for $p \geq 1$
$$V_p = q_{5^{dp}} + k_p.$$
 
 Define then the sequence  ${\lbrace r_n \rbrace}_{n \in \NN} $, successively on the intervals $[V_p,V_{p+1}-1]$, $p=1,\cdots$,  in the following way
\begin{itemize}
\item If $n \in [V_p,V_{p+1}-1]$
for some  $p \geq1$ and is such that $n= q_l+  k_p$, for some $l \in [{5}^{dp},{5}^{d(p+1)}]$,  then take $r_n = l^{-{1/d}},$
\item $r_n = 0, $   otherwise.
\end{itemize}

Then, given any $x_0 \in \TT^d$, take $\B$ to be the sequence of balls centered at $x_0$ with radii $r_n$. Clearly  $\sum r_n^d \geq \sum 1/l = \infty,$ so $\B$ satisfies (\ref{div}).

But it follows from our definition of the sequence $q_l$ that the indexes $n \in [V_p,V_{p+1}-1]$ of the nonempty balls of $\B$, that is  $n=q_l +k_p$ for some $l \in [{5}^{dp},{5}^{d(p+1)}]$, satisfy 
$$\| -n\a_1 +k_p\a_1 \| \leq {e}^{-5^{dp}},$$
together with $r_n = l^{-{1/d}} \leq 5^{-p}$, so that for every $n \in [V_p,V_{p+1}-1]$ we have
$$ T_\a^{-n} B(x_0,r_n) \subset [x_0 -k_p \a_1 - 4^{-p}, x_0 -k_p \a_1 + 4^{-p}] \times \TT^{d-1}.$$
Since the right hand sets  are disjoint for $p \geq 1$ we conclude that \newline 
$\lim \sup T_\a^{-n}  B(x_0,r_n)$ is empty. \carre


\subsection{} \label{24} We now prove that a translation $T_\a$ that is not of constant type does not have the MSTP:
Let $Q_n$ be a  sequence of integers such that for every $n \in \NN$, $Q_{n+1} \geq 2Q_{n}$ and 
\begin{eqnarray} \label{qn} 
{\|Q_n\a\|}_{\ZZ} \leq {1 \over n^{2d+3} Q_n^{1 \over d}}.
\end{eqnarray}

Let $U_n = n^{2d}Q_n$ and $R_n = n^{-2} Q_{n}^{-{1 \over d}}$ and define 
$$r_l = R_n,  \quad  {\rm for \ any \ } l \in [U_{n-1},U_n-1].$$
Clearly $r_l$ is decreasing, while $\sum_{l=U_{n-1}}^{U_n-1}   r_l^d \geq 1/2 U_n R_n^d =1/2$, which implies implies $\sum r_l^d = + \infty$.
 
Next, consider the set
\begin{eqnarray*}   \bigcup_{l=U_{n-1}}^{U_n-1}
  T_\a^{-l} B(x_0, r_l)  &=& \bigcup_{l=U_{n-1}}^{U_n-1}  T_\a^{-l} B(x_0,  R_n) 
 \nonumber \\ 
&\subset&   \bigcup_{l=  0 }^{U_n-1}
 T_\a^{-l} B(x_0,  R_n).  
 \end{eqnarray*}
From (\ref{qn}) we have for any
 $k Q_n$, such that $k Q_n \leq n^{2d}Q_n$
$$ {\| kQ_n\a \|}_{\ZZ} \leq { 1 \over
n^3 Q_n^{{1 \over d}}},$$
consequently, for any $n \geq 1$,
$$\bigcup_{l=  0 }^{U_n-1}
 T_\a^{-l} B(x_0,  R_n) \subset
\bigcup_{l=  0 }^{Q_n-1}
 T_\a^{-l} B(x_0,  2R_n).$$
But 
$$\mu \left( \bigcup_{l=  0 }^{Q_n-1}
 T_\a^{-l} B(x_0,  2R_n)  \right) \leq {C \over n^{2d}}$$
for some constant $C(d)$.

Finally, we deduce that for any $p \geq n \geq 1$
\begin{eqnarray*} 
 \mu \left(\bigcup_{l=U{n-1}}^{U_p-1} T_\a^{-l} B(x_0, r_l)  \right) &\leq& \sum_{k=n}^{p} {C \over n^{2d}} \\
&\leq
& {C' \over n^{2d-1}}, \end{eqnarray*}
for some constant $C'>0$, hence the sequence $\B = {\lbrace B(x_0,r_l) \rbrace}_{l \in \NN}$ is not BC and the second assertion of theorem \ref{theorem1} follows.  \carre

\subsection{} We assume now that $\a \in \RR^d$ is a vector of constant type. Then there exists $\varepsilon(\a)$ such that for any $x_0 \in \TT^d$, we have for all $Q \in \NN$:
\begin{eqnarray}  {\rm \ The \ balls \ } T_\a^{-l} B(x_0,{\varepsilon(\a) / Q^{1 \over d}}), { \rm \ }  l\leq 2Q-1 {\rm \ are \ disjoint} 
\label{disjoint}
\end{eqnarray}

Let $x_0 \in \TT^d$ and adopt the notation $x_l = T_\a^{-l}x_0$ for $l \in \NN$. From remark \ref{eta}, we just have  to show that there exists a constant $\eta(\a)>0$, such that for any decreasing sequence $r_l\geq 0$, satisfying $\sum_{l \geq 0} r_l^d = +\infty$,  we have 
\begin{eqnarray} \label{i} \mu\left(\bigcup_{l\geq 0} B(x_l,r_l)\right) \geq \eta.
\end{eqnarray}

The following elementary lemma will be helpful in  showing that under the condition (\ref{disjoint}) we gain some measure when we consider $\mu (\cup_{l\leq 2Q-1} B(x_l, r_l))$ 
as compared to the measure of $\cup_{l \leq Q} B(x_l, r_l)$, at least if the latter is not already far from zero. We denote by $V(d)$ the volume of the unit ball in $\RR^d$.

\begin{lemm} For any $\epsilon>0$, for any $Q \geq 4$, for any numbers $r_0 \geq \dots \geq r_{2Q-1}$, and for any $2Q$-uple of points $(y_0, \dots,y_{2Q-1}) \in {\TT^{d}}^{2Q}$ such that the balls $B(y_0,\epsilon/Q^{1/d}),\dots,B(y_{2Q-1},\epsilon/Q^{1/d})$ are disjoint, we have at least one of the following alternatives:
\begin{itemize} 
\item[(i)] $\displaystyle{\mu \left(\bigcup_{l=0}^{Q-1}B(y_l,r_l) \right)  \geq V(d) {\left({\epsilon \over 10}\right)}^d.}$ 

\item[(ii)] $\displaystyle{\mu \left( \bigcup_{l=0}^{2Q-1}B(y_l,r_l) \right)   \geq \mu \left( \bigcup_{l=0}^{Q-1}B(y_l,r_l) \right)+ {Q \over 2} \mu\left( B(y_0,r_{2Q-1}) \right)}$.
\end{itemize}
\end{lemm}

\noindent {\sl Proof.} If $r_{Q-1} \geq \epsilon / (10 Q^{1/d})$   we have that 
$ \cup_{l=0}^{Q-1}B(y_l,\epsilon / (10 Q^{1/d})) \subset \cup_{l=0}^{Q-1}B(y_l,r_l)$, but the first is a disjoint union and we obtain  (i).  If to the contrary $r_{Q-1} < \epsilon / (10 Q^{1/d})$,  
we have that if $B(y_l,r_l)  \cap B(y_{l'},r_{l'}) \neq \emptyset$ for some pair $0 \leq l \leq Q-1 < l'$ then  a proportion larger than $1/5$ of $B(y_{l'},\epsilon/Q^{1/d})$  has to be included  in $B(y_l,r_l)$ (this is due to the fact that  $r_{l'} < \epsilon / (10 Q^{1/d})$ and $B(y_l,\epsilon/Q^{1/d}) \cap B(y_{l'},\epsilon/Q^{1/d}) = \emptyset$). So, if there are nonempty intersections between $B(y_{l'},r_{l'})$  and $\cup_{l=0}^{Q-1}B(y_l,r_l)$  for more than $[Q/2]-1$ indexes $Q \leq l' \leq 2Q-1$, (i) will follow again from disjointness of the sequence $B(y_{l'},\epsilon/Q^{1/d}), l'=Q,\dots,2Q-1$. If this is not the case we get (ii) by disjointness of the more than $[Q/2]+1$ balls $B(y_{l'},r_{l'})$ that do not intersect $\cup_{l=0}^{Q-1}B(y_l,r_l)$, and by the fact that $r_l$ is decreasing.  \carre

\vspace{0.2cm} 

Going back to the proof of the MSTP, take a decreasing sequence $r_n \geq 0$ satisfying $\sum_{l \geq 0} r_l^d = +\infty$ and let $\eta := V(d) {({\varepsilon(\a) / 10})}^d$. Assume by contradiction that (\ref{i}) does not hold. Let $n \in \NN$ and apply, in view of (\ref{disjoint}), the above lemma to the $2^{n+1}-$uple of points $x_0,\dots,x_{2^{n+1}-1}$. Our assumption that  alternative (i) fails forces that 
 $$\mu \left(\bigcup_{l=0}^{2^{n+1}-1}B(x_l,r_l) \right) \geq \mu \left(\bigcup_{l=0}^{2^{n}-1}B(x_l,r_l) \right) + 2^{n-1} \mu \left(B(x_0,r_{2^{n+1}-1})\right)$$
which implies by recurrence for $n\geq 2$
 $$\mu \left(\bigcup_{l=0}^{2^{n+1}-1}B(x_l,r_l) \right) \geq \sum_{p=2}^{n}  2^{p-1} \mu \left(B(x_0,r_{2^{p+1}-1}) \right)$$
but the divergence of $\sum  \mu(B(x_0,r_l)) \leq \sum 2^p \mu(B(x_0,r_{2^{p}-1}))$ yields then to a contradiction. \carre

\section{Mixing vs STP}

\subsection{} \label{blabla} Many ergodic properties of a dynamical system can be  characterized in terms of Borel Cantelli 
sequences. We restate some of the results  that were surveyed in \cite{kleinbock-chernov} and refer to the bibliography therein: 
\begin{prop} Let $(T,M,\mu)$ be a dynamical system. Then
\begin{itemize}
\item[ (i)] $(T,M,\mu)$ is ergodic iff every constant sequence $A_n
= A$, $\mu(A) > 0 $, is BC;
\item[ (ii)] $(T,M,\mu)$ is weak mixing  iff every sequence $\A$ that only contains finitely many distinct sets,
none of them of measure zero,  is BC;
\item[ (iii)] $(T,M,\mu)$ is lightly mixing iff every sequence
$\A $ that only contains finitely many distinct sets,
possibly of measure zero, and satisfies (\ref{div}) is BC.
\end{itemize}
\end{prop}

If topological restrictions on the target are added (for example centered
balls for the STP or decreasing centered balls for the MSTP) it 
is possible to let the measure of the target $\mu(A_n)$ go to zero. It was proved for example that expanding maps of the
circle \cite{philipp} and Anosov diffeomorphisms on topological spaces \cite{dolgopyat} had the shrinking target property. 

Spectacular applications of dynamical Borel-Cantelli lemmas appear in the study of geodesic flows on homogeneous spaces. It turns out that establishing the BC property for sequences of shrinking cusp neighborhoods that satisfy (\ref{div})  is sufficient to obtain {\it logarithmic laws} for geodesics and  to derive  number theoretical applications in Diophantine approximations (Cf. \cite{sullivan} and \cite{kleinbock-margulis}).  In \cite{meaucourant}, the MSTP (with closed balls for targets) was obtained for geodesic flows on hyperbolic manifolds of finite volume. Except for the latter result, proving the STP often relies heavily on the exponential decay of correlations of the corresponding transformations (Cf. comments after the definition \ref{bc} above).

It is then natural to show by an example that mixing on its own does not force a system to have the STP, or even the MSTP, which can then be put in contrast with the first part of theorem \ref{theorem1}. A yet stronger property than mixing is given by the following 

\begin{defi}[Mixing of all orders] A dynamical system $(T,M,\mu)$ (or flow $(T^t,M,\mu)$) is said to be {\it mixing of  order  $p \geq 2$} if, for any sequence $\displaystyle{ {\lbrace (u_n^{(1)},\cdots, u_n^{(p-1)}) \rbrace}_{n \in \NN}}$, where for $i=1, \cdots,p-1$ the ${\lbrace u_n^{(i)} \rbrace}_{n \in \NN}$ are sequences of integers (or real numbers) such that $\displaystyle{ \lim_{n \rightarrow  \infty} u_n^{(i)} = \infty}$, and for any $p$-uple $(A_0,\cdots,A_{p-1})$ of measurable subsets of $M$, we have 
$$ \lim_{n\rightarrow \infty} \mu \left(T^{-u_n^{(1)} - \cdots - u_n^{(p-1)}} A_{p-1} \cap \cdots \cap T^{-u_n^{(1)}} A_1 \cap A_0 \right) = \mu(A_{p-1}) \cdots \mu(A_0).$$

The general definition of mixing corresponds to mixing of order 2. A system is said to be {\it mixing of all orders} if it is mixing of order $l$ for any $l \geq 2$.
\end{defi}

\subsection{} We now introduce the setting in which our examples will be produced.

\begin{defi} Given $\a \in \RR^d$ and $\phi \in C^1(\TT^{d+1},\RR^*_+)$ we define the reparametrization of $R_{t(\a,1)}$ by $\phi$ to be the flow given by the differential system on $\TT^{d+1}$
$${dx \over dt} = \phi(x) (\a,1),$$
and that we denote by $T^t_{(\a,1),\phi}$. If $T_\a$ is a minimal translation of $\TT^d$, i.e. if the coordinates of $\a$ are independent over $\QQ$ then $T^t_{(\a,1),\phi}$ preserves a unique probability measure $\mu_{\phi}$ that is equivalent to the Haar measure on the torus and is precisely given by the density ${1/\phi}$.

For $k \in \RR$, we denote by $(\TT^{d+1}, T^k_{(\a,1),\phi},\mu_{\phi})$ the system corresponding to the time $k$ map of the reparametrized flow.
\end{defi}

In \cite{mixing} mixing real analytic reparametrizations of linear flows $R_{t(\a,1)}$ on $\TT^3$ were studied and in \cite{s} it was proved that the reparametrizations could be chosen so as to insure that the flows we obtain are mixing of all order. The vectors $\a \in \RR^2$ involved had both coordinates Liouvillean baring arbitrarily fast power-like simultaneous approximations (i.e. for any $j \in \NN$,  ${\|Q \a\|}_{\ZZ} \leq Q^{-j}$ has infinitely many solutions). Mixing for these flows is very slow and we have

\begin{theo} \label{nostp} Let $\phi \in C^1(\TT^{d+1},\RR^*_+)$. If $\a \in \RR^d$  satisfies $\a \notin \Omega_d(2/d)$ then  $(T^1_{(\a,1),\phi},\TT^{d+1},\mu_{\phi})$ does not have the monotone shrinking target property as stated in definition \ref{tttSTP2}.
\end{theo}

Because it preserves a measure with smooth strictly positive density, the map $T^1_{(\a,1),\phi}$ can be conjugated to a map on the torus that preserves the Haar measure and the conjugacy can be taken as regular as $\phi$ is \cite{moser}. Furthermore, it is not hard to see that the conjugated flow will continue to lack the MSTP.  By the constructions in \cite{mixing,s} we thus obtain
\begin{coro} \label{corr} There exist on $\TT^3$ real analytic maps that preserve the Haar measure, are mixing of all orders, and 
do not have the monotone shrinking target property.
\end{coro}
 

\subsection{\sl Proof of theorem \ref{nostp}.} Let $x = (x_1,\cdot,x_{d+1}) \in \TT^{d+1}$ and denote by ${x}^{(d)}$ 
the point $(x_1,\cdots,x_d) \in \TT^d$. The $d$-dimensional torus $\TT^d \times \lbrace x_{d+1} \rbrace$ is a global section of the flow $T^t_{(\a,1),\phi}$  on which the Poincar\'e return map is  the translation $T_\a$ and the return time of any point is bounded between two strictly positive constants $c$ and $C$.\footnote{These constants are the {\it suprema} of the return time function given by \newline $\varphi(x) = \int_{\TT^1} {1 / \phi(x+s\a,s)}ds.$}. Therefore we have for any ball $B(x,r) \subset \TT^{d+1}$ and for any $\tau \geq 0$
\begin{eqnarray} \label{special} 
\bigcup_{t=-\tau}^{0} T^t_{(\a,1),\phi} B(x,r) 
\subset  \bigcup_{s=0}^{C} T^{-s}_{(\a,1),\phi}  \bigcup_{l=0}^{\left[{\tau \over c}\right]} T_\a^{-l} B(x^{(d)},r). 
\end{eqnarray}

Since $\a \notin \Omega_d(2/d)$, let $Q_n$ be a sequence of integers such that $Q_{n+1} \geq 2Q_n$ and 
\begin{eqnarray} \label{qn2} {\|Q_n \a \|}_{\ZZ} \leq  {1 \over n^{2d+5} Q_n^{2 \over d}}. \end{eqnarray}


Denote for $n \geq 1$
\begin{eqnarray*} R_n &=& {1\over n^{2}Q_n^{1 \over d}} \\
U_n &=& {n}^{2d+2}[Q_{n}^{1+{1\over d}}] \end{eqnarray*}
and define a sequence of radii as follows 
$$r_l = R_n,  \quad  {\rm for \ any \ } l \in [U_{n-1},U_n-1].$$


Clearly $r_l$ is decreasing, while $\sum_{l=U_{n-1}}^{U_n-1} r_l^{d+1} \geq 1/2$, which implies
$\sum r_l^{d+1} = + \infty$. The latter implies that the sequence $B_l(x,r_l)$ satisfies (\ref{div}) for $\mu_\phi$ that is equivalent to the Haar measure on $\TT^{d+1}$.

But  (\ref{qn2}) implies, as in \S \ref{24}, that 
\begin{eqnarray} \label{2/d} 
A_n= \bigcup_{l=0}^{[{U_n \over c}]} T_\a^{-l} B(x^{(d)},R_n)  \subset \bigcup_{l=0}^{Q_n}  T_\a^{-l} B(x^{(d)},2R_n). 
\end{eqnarray}
Denoting by $\mu^{(d)}$ the Haar measure on $\TT^{d}$ we have from (\ref{2/d})
$$\mu^{(d)} (A_n) \leq {C_1 \over n^{2d}}$$
for some constant  $C_1 >0$. Hence
$$\mu^{(d+1)} (\bigcup_{s=0}^{C} T^{-s}_{(\a,1),\phi} A_n ) \leq {C_2 \over n^{2d}}$$ 
which also implies by equivalence of the measures $\mu^{(d+1)}$ and $\mu_\phi$ 
$$\mu_{\phi} (  \bigcup_{s=0}^{C} T^{-s}_{(\a,1),\phi} A_n) \leq {C_3 \over n^{2d}}.$$ 
In view of (\ref{special}) we have  
\begin{eqnarray*} \bigcup_{l=-U_n+1}^{-U_{n-1}} T^{l}_{(\a,1),\phi} B(x,R_n)   &\subset& \bigcup_{t=-U_n}^{0}  T^{l}_{(\a,1),\phi} B(x,R_n)  \\
&\subset& \bigcup_{s=0}^{C} T^{-s}_{(\a,1),\phi} A_n \end{eqnarray*}
and conclude that
$$\mu_\phi \left(\bigcup_{l=-U_n+1}^{-U_{n-1}} T^{l}_{(\a,1),\phi} B(x,r_l) \right) \leq {C_3 \over n^{2d}},$$
which implies that $\lim \sup T^{-n}_{(\a,1),\phi} B(x,r_n)$ has zero measure. \carre
  
\begin{rema} The same arguments made here on absence of the MSTP could be carried out for the mixing flows on $\TT^2$ introduced in \cite{Kochergin2} under a suitable restriction of their rotation vector. Mixing for all orders however, is not yet established for these flows. 
\end{rema}


Finally, it is worth noticing while we compare the positive result related to constant type translations in theorem \ref{theorem1} and the counterexamples of corollary \ref{corr}, that  the monotone shrinking target property and the mixing property, that can both be interpreted as a qualitative strengthening of ergodicity, are actually independent.

\frenchspacing
\bibliographystyle{plain}
 
\end{document}